\documentclass[12pt]{article}

\usepackage{amssymb}

\newcommand{\be}{\begin{equation}}
\newcommand{\ee}{\end{equation}}
\newcommand{\bea}{\begin{eqnarray*}}
\newcommand{\eea}{\end{eqnarray*}}
\newcommand{\ba}{\begin{array}}
\newcommand{\ea}{\end{array}}
\newcommand{\bi}{\begin{itemize}}
\newcommand{\ei}{\end{itemize}}
\newcommand{\bc}{\begin{center}}
\newcommand{\ec}{\end{center}}
\newcommand{\bfr}{\begin{flushright}}
\newcommand{\efr}{\end{flushright}}
\newcommand{\f}{\frac}
\newcommand{\ov}{\overline}

\newcommand{\ds}{\displaystyle}

\newcommand{\q}{\quad}

\newcommand{\bb}{\mathbb}

\begin{document}

\title{On An Identity Derived From Interpolation Theory}
\author{Sorin G. Gal\\
Department of Mathematics,\\
University of Oradea,\\
3700 Oradea, Romania\\
e-mail: galso@math.uoradea.ro}
\date{}
\maketitle

The aim of this note is to show how can be derived from the properties
of fundamental interpolation polynomials some identities.
Firstly let us recall some known facts in interpolation theory.
Let $f:[a,b]\to\bb{R}$ and $x_1<\dots<x_n$ be distinct points
(knots) in $[a,b]$.
It is well-known that the Hermite-Fej\'er interpolation polynomial
$H_{2n-1}(f)(x)$ (of degree $2n-1$) attached to $f$ on the knots
$x_i,\ i=\ov{1,n}$, satisfies
$H_{2n-1}(f)(x_i)=f(x_i)$,
$H'_{2n-1}(f)(x_i)=0$, $i=\ov{1,n}$,
$H_{2n-1}(f)(x)=\ds\sum_{i=1}^n h_{i,n}(x)f(x_i)$,
where
$\ds\sum_{i=1}^n h_{i,n}(x)=1$, for any $x\in\bb{R}$.
It follows
\be
\sum_{i=1}^n h_{i,n}^{(p)}(y_0)=0,\ \forall \ p\in\bb{N},\ y_0\in\bb{R}.
\ee

The idea is that by using (1) for some special choices for
$x_i$, $i=\ov{1,n}$, $y_0$ and $p\in\bb{N}$,
to get some interesting identities.

In this sense let us present the following.

{\bf Application.}
Prove that for any odd number $n\ge 3$, the identity
\be
\sum_{k=1}^\f{n-1}{2} \f{2}{\sin^2\ds\f{k\pi}{n}}=\f{n^2-1}{3},
\ee
holds.

{\bf Proof.}
We use (1) with $p=2,\ y_0=0$ and $x_i=\cos\ds\f{2i-1}{2n}\pi,\
i=\ov{1,n}$ (the Chebyshev knots of first kind).
It is known that
$$h_{i,n}(x)=\f{1}{n^2}\left[\f{T_n(x)}{x-x_i}\right]^2 (1-xx_i),\q
i=\ov{1,n},$$
where $T_n(x)=\cos[n\arccos x]$ are the Chebyshev polynomials
of first kind.

Let $n\ge 3$ be odd.
After some simple calculations (we use here some known results
on $T_n(x)$ in e.g. [1,p.213-214])
$$h''_{i,n}(0)=\f{2}{n^2}\left[\f{T'_n(0)}{x_i}\right]^2=\f{2}{x_i^2},
\q\mbox{if}\q i\ne\f{n+1}{2}$$
and
$$h''_{\f{n+1}{2},n}(0)=\f{2}{3}(1-n^2).$$

Since
$h''_{i,n}(0)=h''_{n+1-i}(0)$,
$i\in\{1,2,\dots,n\}\setminus\left\{\ds\f{n+1}{2}\right\}$,
applying (1) we immediately obtain (2).

{\bf Remarks.}
1) Try a direct proof for (2).

2) Other identities can be derived by taking in (1) the same knots
$x_i$ as above, $y_0=0$, but different values for $p>2$.
This challenge is left to the reader.

3) Another open question would be to derive from (1) similar
identities for more general knots $x_i,\ i=\ov{1,n}$,
as for example the Jacobi knots.


\begin{thebibliography}{99}
\bibitem{1}
Gh. Mocic\u a,
{\it Problems of special functions} (in Romanian),
Ed. Did. Ped., Bucharest, 1988.
\end{thebibliography}
\end{document}